\documentclass[12pt,noamsfonts]{amsart}

\usepackage{amsmath}
\usepackage{amsxtra}
\usepackage{amsthm}
\usepackage{amsfonts}
\usepackage{amssymb}
\usepackage{amscd}
\usepackage{eucal}

\newcommand{\marginlabel}[1]%
  {\mbox{}\marginpar{\raggedleft\hspace{0pt}\bfseries\sf#1}}

\def\ZZ{{\mathbb Z}}
\def\N{{\mathbb N}}

\def\CC{{\mathbb C}}
\def\AA{{\mathbb A}}

\def\QQ{{\mathbb Q}}
\def\PP{{\mathbb P}}

\def\cI{\mathcal{I}}

\def\cO{\mathcal{O}}
\def\cM{\mathcal{M}}

\DeclareMathOperator{\Hom}{Hom}

\DeclareMathOperator{\Spec}{Spec}

\newtheorem{lemma}{Lemma}[section]
\newtheorem{theorem}[lemma]{Theorem}
\newtheorem{corollary}[lemma]{Corollary}
\newtheorem{proposition}[lemma]{Proposition}

\theoremstyle{definition}

\newtheorem{example}[lemma]{Example}
\newtheorem{question}[lemma]{Question}
\newtheorem{conjecture}[lemma]{Conjecture}
\newtheorem{remark}[lemma]{Remark}

\numberwithin{equation}{section}

\newcommand{\bean}{\begin{eqnarray}}
\newcommand{\eean}{\end{eqnarray}}
\newcommand{\be}{\begin{displaymath}}
\newcommand{\ee}{\end{displaymath}}
\newcommand{\bea}{\begin{eqnarray*}}   
\newcommand{\eea}{\end{eqnarray*}}

\newcommand{\thmref}[1]{Theorem~\ref{#1}}

\newcommand{\propref}[1]{Proposition~\ref{#1}}

\newcommand{\nc}{\newcommand}
\nc{\on}{\operatorname}
\nc{\Z}{{\mathbb Z}}
\nc{\C}{{\mathbb C}}
\nc{\pa}{D}
\nc{\F}{{\mathcal F}}
\nc{\arr}{\rightarrow}
\nc{\larr}{\longrightarrow}
\nc{\la}{\lambda}
\nc{\g}{{\mathfrak g}}
\nc{\Ga}{\Gamma}
\nc{\wt}{\widetilde}
\nc{\wh}{\widehat}
\nc{\el}{\ell}
\nc{\bi}{\bibitem}
\nc{\om}{\omega}
\nc{\ol}{\overline}
\nc{\ds}{\displaystyle}
\nc{\GL}{^L\G}
\nc{\gL}{^L\g}
\nc{\mc}{\mathcal}
\nc{\Spe}{\on{Spec}}
\nc{\NN}{\mc N}
\nc{\hh}{{\mathfrak h}}
\nc{\K}{{\mc K}}
\nc{\sw}{{\mathfrak s}{\mathfrak l}}

\begin{document}

\begin{center}

{\bf Jet Schemes of Locally Complete Intersection}

{\bf Canonical Singularities} 

{ by \bf Mircea Musta\c{t}\v{a}}

\smallskip

{ with an appendix}

{by {\bf David Eisenbud} and {\bf Edward Frenkel}}

\date{\today}
\end{center}

\section*{Introduction}

Let $X$ be a variety defined over an algebraically closed field $k$
of characteristic zero. The $m$th jet scheme $X_m$ of $X$ is a scheme
 whose closed points over $x\in X$ are morphisms
$\cO_{X,x}\longrightarrow k[t]/(t^{m+1})$.
 When $X$ is a smooth variety,
this is an affine bundle over $X$, of dimension $(m+1)\dim\,X$.
The space of arcs $X_{\infty}$ of $X$ is the projective limit
$X_{\infty}=\projlim_mX_m$.

Our main result is a proof of the following theorem, which was conjectured by
Eisenbud and Frenkel:

\begin{theorem}\label{main_conj}
If $X$ is locally a complete intersection variety, then $X_m$ is
irreducible for all $m\geq 1$ if and only if $X$ has rational singularities.
\end{theorem}

In the appendix, Eisenbud and Frenkel apply this result when $X$
is the nilpotent cone of a simple Lie algebra to extend results of
Kostant in the setting of jet schemes.

Note that since
 $X$ is assumed to be locally complete intersection, hence Gorenstein,
a result of Elkik \cite{elkik} and Flenner \cite{flenner} says that $X$
has rational singularities if and only if it has canonical singularities.

We make also make the following conjecture toward a similar characterization
of log canonical singularities.

\begin{conjecture}\label{main_conj2}
If $X$ is locally a complete intersection, normal variety, then
$X_m$ is pure dimensional for all $m\geq 1$ if and only if $X$ has log
canonical singularities.
\end{conjecture}

We prove the ``only if'' part of Conjecture~\ref{main_conj2} and show that
the ``if'' part is equivalent to a special case
of the Inverse of Adjunction Conjecture 
due to Shokurov and several other people (see \cite{kollar}, Conjecture 7.3). 

\smallskip

One should contrast Theorem~\ref{main_conj} with the following result
of Kolchin.

\begin{theorem}[\cite{kolchin}]
If $X$ is a variety over a field of characteristic zero, then
$X_{\infty}$ is irreducible.
\end{theorem}

However, when $X$ is locally a complete intersection and has rational
singularities, Theorem~\ref{main_conj} gives much more information about
$X_{\infty}$ (for example, as we will see, it implies that $X_{\infty}$
is reduced). 

The main technique we use in proving Theorem~\ref{main_conj} is 
motivic integration, as developed by Kontsevich, Denef and Loeser,
and Batyrev.
Here is a brief description of 
the proof of Theorem~\ref{main_conj}. Consider
an embedding $X\subset Y$, of codimension $r$,
 where $Y$ is smooth, and an embedded resolution
of singularities $\gamma\,:\,\widetilde{Y}\longrightarrow Y$ for $X$.
There is a function $F_X$ on $Y_{\infty}$, defined by
$F_X(w)={\rm ord}(w(\cI_{X,y}))$, where $w$ is considered 
as a morphism  $w\,:\,\cO_{Y,y}\longrightarrow k[[t]]$.
 By integrating the function
$f\circ F_X$ on $Y_{\infty}$, for a convenient function $f\,:\,\N
\longrightarrow\N$, we get a Laurent series in two variables which
encodes information about the dimensions of $X_m$ and the number of irreducible
components of maximal dimension. 

Applying the change of variable formula in \cite{batyrev} or \cite{denef1},
this integral can be expressed as an integral on $\widetilde{Y}_{\infty}$
and since $\gamma^{-1}(X)$ is a divisor with normal crossings, this can
be explicitely computed. If $\gamma^{-1}(X)=\sum_{i=1}^ta_iE_i$, where
$E_1$ is the only exceptional divisor dominating $X$ and the
discrepancy of $\gamma$ is $W=\sum_{i=1}^tb_iE_i$, then we see that
$b_j\geq ra_j$ for all $j\geq 2$ if and only if $\dim\,X_m=(m+1)\dim\,X$,
and $X_m$ has exactly one component of maximal dimension for all $m$.
When $X$ is locally a complete intersection, this says precisely that $X_m$ is
irreducible for all $m$.

The last step needed is that this numerical condition is equivalent with
$X$ having
canonical singularities when $X$ is locally a complete intersection. We 
consider the following construction of $\gamma$: let 
$p\,:\,B\longrightarrow Y$
 be the blowing-up of $Y$ along $X$, $F$ the exceptional   
divisor, and $\tilde{p}\,:\,\widetilde{Y}\longrightarrow B$
an embedded resolution of singularities for $F\subset B$. We take
$\gamma=p\circ\tilde{p}$ and we show that the numerical condition is
equivalent with $(B,F)$ being canonical. By a result of Stevens \cite{stevens},
this is equivalent with $F$, hence $X$ having canonical singularities.

The computation of motivic integrals gives an analogous
condition  
for a variety which is locally a complete intersection
 to have pure dimensional jet schemes. The condition is that 
$b_j\geq ra_j-1$, for all $j$. Conjecture~\ref{main_conj2}
can therefore be translated into a conjectural analogue
of the result of Stevens for log canonical singularities.

The technique we use to describe singularities in terms of jet schemes
can be applied also to study pairs $(X,D)$, where $X$ is a smooth variety
and $D$ an effective $\QQ$-divisor on $X$. For example, we prove in
\cite{mustata} the following characterization of log canonical pairs.

\begin{theorem}
Let $X$ be a smooth variety and $D$ an effective divisor on $X$ with
integral coefficients.
\item{i)} For every positive integer $n$, the pair $(X,{1\over n}D)$
is log canonical if and only if
$$\dim\,D_m\leq (m+1)(\dim\,X- 1/n),$$
for all $m$.
\item{ii)} The log canonical threshold of $(X,D)$ is given by
$$c(X,D)=\dim\,X-\sup_{m\geq 0}{{\dim\,D_m}\over {m+1}}.$$
\end{theorem}

\smallskip

In the first section of the paper we give the definition of jet schemes
and discuss the irreducibility condition for jets of locally complete
intersection varieties. The condition that the jet schemes are irreducible
(or pure dimensional) can be formulated in terms of the dimension of
the space of jets lying over the singular part.
 
Starting with equations of $X$ in an affine space $\AA^N$, it is easy to
give equations for $X_m$. If $X$ is defined by $(f_{\alpha}(U))
\subset k[U_1,\ldots,U_N]$, then $X_m$ can be defined in $\AA^{(m+1)N}$
by $(f_{\alpha}^{(j)})_{0\leq j\leq m}$,
where $f_{\alpha}^{(j)}\in k[U,U',\ldots,U^{(m)}]$ is given by
$f_{\alpha}^{(j)}=D^j(f_{\alpha})$, $D$ being the derivation 
taking $U_i^{(l)}$ to $U_i^{(l+1)}$ for $l\leq m-1$. 

Using this explicit description and the irreducibility criterion, one can
check directly that some jet schemes are irreducible. We give applications
in the last section.

We deduce from
the description by equations
 that if $X$ is locally a complete intersection and $X_m$ is
pure dimensional, then it is locally a complete intersection, too. It follows
that if $X_m$ is irreducible, then it is reduced. We see also
that if $X_m$ is irreducible, then so is $X_{m-1}$.

In the second section we show that the numerical condition coming
from an embedded resolution of singularities of $X$ is equivalent
with $X$ having canonical singularities. The third section
uses motivic integration as we described.

In the last section we discuss several examples and open questions.
We consider first the small dimension case. If $X$ is a singular curve,
then $X_m$ is reducible for all $m$. If $X$ is a surface (and
${\rm char}\,k=0$), we show that being locally a complete intersection
is a necessary condition for the irreducibility of jet schemes.
More precisely, if $X$ is a surface, then $X_m$ is irreducible
for all $m$ if and only if all the singular points of $X$
are rational double points.

We give an example of a toric variety of dimension $3$ which shows
that in Theorem~\ref{main_conj} it is not possible to replace 
locally complete intersection with Gorenstein. On the other hand,
the example of the cone over the Segre embedding of $\PP^1
\times\PP^n$, with $n\geq 2$, shows that the condition of being
locally a complete intersection is not necessary in order to have all
the jet schemes irreducible.

Our results in the second and the third section, where we used the theory
of singularities of pairs and motivic integration, rely on the fact that
the characteristic of the ground field is zero.
We discuss briefly a possible analogue of Theorem~\ref{main_conj}
in positive characteristic and we 
 end with a characteristic free proof
of the fact that if $X$ is a locally complete intersection toric variety,
then $X_m$ is irreducible for all $m$. This is achieved
 using an inductive
description due to Nakajima \cite{nakajima} for such varieties
in order to describe a desingularization of the ``dual'' toric variety.

\subsection{Acknowledgements}
It is a pleasure to thank David Eisenbud and Edward Frenkel who got
me interested in this problem and helped me all along with useful
suggestions and comments. Without their constant encouragement and support,
this work would have not been done. I am grateful to 
Joe Harris, J\'{a}nos Koll\'{a}r,
Monique Lejeune-Jalabert, Miles Reid
 and Matthew Szczesny for their help and suggestions
during various stages of this project. Thanks are due also to Lawrence Ein
who pointed out some incomplete arguments in an earlier version of this paper.

\section{Jet basics}

The study of singularities via the space of arcs has gotten a lot of attention 
recently. Nash initiated this study in \cite{nash}. He
suggested that the study of the images $\eta_m(X_{\infty})\subseteq X_m$
for all $m$, where $\eta_m$ are the canonical projections, should give
information about the fibers over the singular points in the 
desingularizations of $X$.
For more on this approach, see \cite{lejeune} and \cite{lejeune-reguera}.
For applications of spaces of arcs with a different flavour, 
for example the
proof of a geometric analogue of Lang's Conjecture, see \cite{buium}.

\smallskip

We start by reviewing the definition and the general properties of 
jet schemes. In the case of locally complete intersection varieties, we
give an irreducibility criterion for these schemes and show that under 
the irreducibility assumption, they are, as well, locally complete
intersection varieties.

Let $k$ be an algebraically closed field.
If ${\mathcal Sch}/k$ is the category of schemes of finite type over $k$, 
consider for every $m\geq 0$ the covariant functor 
${\bf F}\,:{\mathcal Sch}/k\longrightarrow {\mathcal Sch}/k$, given by
${\bf F}(Y)=Y\times_{\Spec\, k}{\Spec\,k[t]/(t^{m+1})}$.
This functor has a right adjoint ${\bf G}\,:\,{\mathcal Sch}/k
\longrightarrow {\mathcal Sch}/k$, given by ${\bf G}(X)=X_m$, and
$X_m$ is called the scheme of jets of order $m$ of $X$.

 For an affine scheme
$Y=\Spec\,A$, the adjointness relation says that the $A$-valued points
of $X_m$ are in bijection with the $A[t]/(t^{m+1})$-valued points of $X$.
In particular, there are canonical isomorphisms $X_0\simeq X$
and $X_1\simeq T\,X$, where $T\,X$ is the total tangent space of $X$.

By adjointness, the canonical embeddings $$Y\times_{\Spec \,k}\Spec\,
k[t]/(t^m)\hookrightarrow Y\times_{\Spec \,k}\Spec\,k[t]/(t^{m+1})$$
induce canonical morphisms $\phi_m^X\,:\,X_m\longrightarrow X_{m-1}$,
for $m\geq 1$. We will use the notation
$\pi^X_m=\phi^X_1\circ\ldots\circ\phi^X_m\,:X_m\longrightarrow X$,
but we will supress the variety $X$, whenever this leads to no confusion.

The space of arcs of $X$, denoted by $X_{\infty}$, is the inverse
limit of $\{X_m\}_{m\geq 0}$. This is a scheme over $k$, in general
not of finite type, whose $A$-valued points are in natural bijection
with the $A[[t]]$-valued points of $X$. 

\begin{proposition}\label{etal_com}
If $f\,:\,X\longrightarrow Y$ is an \'{e}tale morphism, then
$X_m\simeq Y_m\times_YX$, for all $m$.
\end{proposition}

\begin{proof}
The assertion follows by adjointness from the fact that
$f$ is also formally \'{e}tale.
\end{proof}

In particular, the construction of jet schemes is compatible with
open immersions. Therefore, in order to describe $X_m$, we may
restrict ourselves to the affine case: suppose $X\subseteq\AA^N$,
$X=\Spec(R)$ and $R=k[U_1,\ldots,U_N]/(f_1,\ldots,f_r)$.
For every $A$, an element in \break
 $\Hom(\Spec\,A[t]/(t^{m+1}),X)$
is given by a morphism $\theta\,:\,k[U_1,\ldots,U_N]
\longrightarrow A[t]/(t^{m+1})$ such that $\theta(f_{\alpha})=0$
for all $\alpha$. The condition $\theta(f_{\alpha})=0$ is equivalent
with the vanishing of each of the coefficients of $t^i$ in
 $\theta(f_{\alpha})$, for $0\leq i\leq m$ and gives $m+1$ equations.
Therefore the map $\theta\longrightarrow
(\theta^{(j)}_i)$, where $\theta(U_i)=\sum_{j=0}^m
\theta^{(j)}_it^j$,
 induces a closed immersion $X_m\hookrightarrow\AA^{(m+1)N}$,
 such that $X_m$ can be defined
by $(m+1)r$ equations.

When ${\rm char}\,k>m$, by normalizing the variables,
 the equations defining $X_m$ can be
written as follows. Let $S_m=k[U_i^{(j)};1\leq i\leq N,
0\leq j\leq m]$ be the coordinate ring of $\AA^{(m+1)N}$
and $D\,:\,S_m\longrightarrow S_{m+1}$ the unique derivation over $k$
such that $D(U^{(j)}_i)=U^{(j+1)}_i$ for all $i$ and $j$.
If we embed $X_m$ in $\AA^{(m+1)N}$ by $\theta\longrightarrow
(j!\,\theta_i^{(j)})$, then
the ideal of $X_m$ is generated by
$f^{(j)}_{\alpha}$, for all $\alpha$ and all $j$, $0\leq j\leq m$,
where $f^{(j)}=D^j(f)$.

We will need later the following lemma.

\begin{lemma}\label{fibers}
For every scheme $X$ and every $u\in X_m$,
either $\phi_{m+1}^{-1}(u)=\emptyset$, or $\phi_{m+1}
^{-1}(u)\simeq\phi_1^{-1}(x)$, where $x=\pi_m(u)$.
\end{lemma}

\begin{proof}
If we look at $A$-valued points, then
${\rm Hom}({\rm Spec}\,A, \phi_{m+1}^{-1}(u))$ is a pseudotorsor
over ${\rm Hom}({\rm Spec}\,A,T_xX)$. Indeed, suppose that
$v\,:\,\cO_{X,x}\longrightarrow A[t]/(t^{m+2})$, $v(y)
=\sum_{i\leq m+1}v_i(y)t^i$ corresponds to an $A$-valued point
of $\phi_{m+1}^{-1}(u)$.  Any other such morphism is of the form
$v'(y)=v(y)+w(y)t^{m+1}$, where $v_0+wt$ is an $A$-valued point of
$T_xX$.

If $\phi_{m+1}^{-1}(u)\neq\emptyset$, then for a fixed closed point
in $\phi_{m+1}^{-1}(u)$, we get induced $A$-valued points in $\phi_{m+1}
^{-1}(u)$, and therefore an isomorphism $\phi_{m+1}^{-1}(u)\simeq T_xX$.
\end{proof}

It is well-known that if $X$ is a smooth, connected variety
of dimension $n$, then  for every $m$, the morphism $\pi_m$ is an 
affine bundle with fiber $\AA^{mn}$. Under these circumstances, $X_m$
is smooth, connected, of dimension $(m+1)n$.

\smallskip

We define now a morphism $\Psi_m\,:\,\AA^1\times X_m\longrightarrow
X_m$ by defining it on $A$-valued points. If $(a, f)$ corresponds to an
$A$-valued point of $\AA^1\times X_m$, where $a\in A$ and
$f\,:\,\Spec\,A[t]/(t^{m+1})\longrightarrow X$, then $\Psi_m(a,f)$
is given by the composition 
$$\Spec\,A[t]/(t^{m+1})\overset{g_a}{\longrightarrow}
\Spec\,A[t]/(t^{m+1})\overset{f}{\longrightarrow} X,$$
where $g_a$ corresponds to the morphism of $A$-algebras which maps $t$ to $at$.

It is clear that these morphisms are compatible with the projections
$\phi_m\,:\,X_m\longrightarrow X_{m-1}$ and that the restriction to
$k^*\times X_m$ defines an action of $k^*$ on $X_m$. Notice that
there is a canonical section of the projection $\pi_m$. This is
the morphism $s_m\,:\,X\longrightarrow X_m$, which takes an $A$-valued
point of $X$, $f\,:\,\Spec\,A\longrightarrow X$
to the composition $f\circ p_m$, where $p_m$ corresponds to the inclusion
$A\hookrightarrow A[t]/(t^{m+1})$. It follows from the definition
that $\phi_{m+1}\circ s_{m+1}=s_m$ and $\Psi_m\vert_{\{0\}\times X_m}
=s_m\circ\pi_m$.

We want to study the irreducible components of $X_m$.
The following  lemma allows us to relate the irreducible components of $X_m$
and $X_{m+1}$.

\begin{lemma}\label{irred_comp}
If $Z\subset X$ is a closed subscheme and $S$ is an
irreducible component of $\pi_m^{-1}(Z)$, then
$\overline{s_m(\pi_m(S))}\subseteq S$. In particular,
$\phi_{m+1}^{-1}(S)\ne\emptyset$.
\end{lemma}

\begin{proof}
It is enough to notice that since $S$ is an irreducible component of
$\pi_m^{-1}(Z)$, we have $\Psi_m(\AA^1\times S)\subseteq S$
and therefore $s_m(\pi_m(S))\subseteq S$.
 \end{proof}

{}From now on, we will restrict ourselves to the case when
$X$ is locally a complete intersection 
(l.c.i. for short) variety. As usual, a variety is
an integral scheme of finite type over $k$. We
denote the smooth part of $X$ by $X_{\rm reg}$ and its complement
by $X_{\rm sing}$.

\begin{proposition}\label{irred_crit}
Let $X$ be an l.c.i. variety of dimension $n$
and $m$ a positive integer. The scheme $X_m$ is pure dimensional
if and only if $\dim\,X_m\leq n(m+1)$, and in this case $X_m$ is
locally complete intersection, too. Similarly, $X_m$ is irreducible
if and only if $\dim\,\pi_m^{-1}(X_{\rm sing})<n(m+1)$.
\end{proposition}

\begin{proof}
We have a decomposition
$$X_m=\pi_m^{-1}(X_{\rm sing})\cup\overline{\pi_m^{-1}(X_{\rm reg})}$$
and in general $\overline{\pi_m^{-1}(X_{\rm reg})}$ is an irreducible
component of $X_m$ of dimension $n(m+1)$. Therefore the ``only if''
part of both assertions is obvious and holds without the l.c.i.
hypothesis.

Suppose now that $\dim\,X_m\leq n(m+1)$. Working locally, we may assume that
$X\subset\AA^N$ and that $X$ is defined by $N-n$ equations. We have seen
that $X_m\subset\AA^{N(m+1)}$ is defined by $(N-n)(m+1)$ equations,
and therefore every irreducible component of $X_m$ has dimension
at least $n(m+1)$. We deduce that $X_m$ is pure dimensional and locally
a complete intersection.

If $\dim\,\pi_m^{-1}(X_{\rm sing})<n(m+1)$, this implies that
$\dim\,X_m\leq n(m+1)$, so that $X_m$ is pure dimensional. The above 
decomposition of $X_m$ shows that $X_m$ is irreducible.
\end{proof}

\begin{proposition}\label{reduced}
If $X$ is an l.c.i.  variety and 
$X_m$ is irreducible for some $m\geq 1$, then $X_m$ is also reduced.
\end{proposition}

\begin{proof}
By Proposition~\ref{irred_crit}, $X_m$ is l.c.i.,
hence Cohen-Macaulay. Since $\pi_m^{-1}(X_{\rm reg})$ is smooth,
$X_m$ is generically reduced, and we conclude by Macaulay's
theorem (see \cite{eisenbud}, Corollary 18.14).
\end{proof}

\begin{proposition}\label{one_implies_lower}
If $X$ is an l.c.i. variety of dimension $n$ and $Z\subseteq X$
is a closed subscheme, then $\dim\,(\pi_{m+1}^{-1}(Z))
\geq\dim\,(\pi_m^{-1}(Z))+n$, for every $m\geq 1$. 
 In particular,
if  $X_{m+1}$ is irreducible or pure dimensional, then so is $X_m$.
\end{proposition}

\begin{proof}
Again, we may assume that $X\subseteq\AA^N$ is defined by $N-n$ equations.
It follows from the equations of $X_{m+1}$ that we have
$\pi_{m+1}^{-1}(Z)\hookrightarrow
\pi_m^{-1}(Z)\times\AA^N$, such that $\pi_{m+1}^{-1}(Z)$ is defined by
$N-n$ equations. The first assertion follows from this once we notice
that by Lemma~\ref{irred_comp}, for every irreducible component
$S$ of $\pi_m^{-1}(Z)$, we have $\phi_{m+1}^{-1}(S)\ne\emptyset$.
The last statement is a consequence of Proposition~\ref{irred_crit}.
\end{proof}

\begin{proposition}\label{normal}
If $X$ is an l.c.i. variety and 
$X_m$ is irreducible for some $m\geq 1$, then $X$ is normal.
\end{proposition}

\begin{proof}
Since $X$ is in particular Cohen-Macaulay, by Serre's Criterion
(see \cite{eisenbud}, Theorem 11.5) it is enough to show that
${\rm codim}(X_{\rm sing},X)\geq 2$. If $X_m$ is irreducible, by
Proposition~\ref{one_implies_lower}, we may assume that $m=1$.
But if ${\rm codim}(X_{\rm sing},X)=1$, since 
for every $x\in X_{\rm sing}$ we have $\dim\pi_1^{-1}(x)=
\dim\,T_xX\geq\dim\,X+1$, it follows
that $\dim\pi_1^{-1}(X_{\rm sing})\geq 2\dim\,X$, contradicting
 Proposition~\ref{irred_crit}.
\end{proof}

\section{A Criterion for L.C.I. Varieties to Have Canonical
Singularities}

In this section we establish the criterion we will use to check
that an l.c.i. variety $X$ has canonical singularities. We embed $X$
in a smooth variety $Y$ and our criterion is in terms of the data
coming from an embedded resolution of singularities of $X\subseteq Y$.

We assume that the characteristic of the ground field is zero.
For the definitions of singularities of pairs, we refer to
\cite{kollar} or to \cite{komo}, Chapter 2.3.
Let $X$ be a normal l.c.i. variety and we fix 
an arbitrary embedding $X\hookrightarrow Y$,
 where $Y$ is a smooth variety.
Let $r$ be the codimension of $X$ in $Y$.

Consider the blowing-up $p\,:\,B=Bl_XY\longrightarrow Y$ 
of $Y$ along $X$, and let $F=p^{-1}(X)$ be the exceptional divisor.
Since $X$ is locally a complete intersection, $F$
is a projective bundle over $X$. 
In particular, $F$ is an integral divisor on $B$, and is locally
a complete intersection. Moreover, $F$ is normal since $X$ is, and 
therefore $B$ is normal, too.

By Hironaka's embedded resolution of singularities
(see \cite{hironaka}), there is a morphism 
$\tilde p\,:\,\widetilde{Y}\longrightarrow B$ which is proper,
an isomorphism over the complement 
of a proper closed subset of $F$, and such that
$\widetilde{Y}$ is smooth and ${\tilde p}^{-1}(F)$ is a divisor with
normal crossings.

Let $\gamma$ be the composition $p\circ\tilde{p}$.
We can write $\gamma^{-1}(X)={\tilde p}^{-1}(F)=E_1+\sum_{i=2}^ta_iE_i$,
where $E_2,\ldots,E_t$ are the exceptional divisors of $\tilde p$,
and $E_1$ is the proper transform of $F$. 

The discrepancy $W$ of $\gamma$
 is defined by the formula $K_{\widetilde{Y}}
=\gamma^{-1}(K_Y)+W$. We write $W=\sum_{i=1}^tb_iE_i$. The following
is our criterion for $X$ to have canonical singularities. 

\begin{theorem}\label{num_can}
With the above notation, $X$ has canonical singularities if and only if
$b_i\geq ra_i$ for every $i\geq 2$.
\end{theorem}

\begin{proof}
Notice first that we have $K_B=p^{-1}(K_Y)+(r-1)F$. Indeed, in order to
compute the coefficient of $F$, we may restrict to an open subset whose
intersection with $X$ is nonempty and smooth, in which case the formula
is well-known.

Consider the divisor $R$
on $\widetilde{Y}$, defined by
$K_{\widetilde{Y}}=\tilde{p}^{-1}(K_B+F)+R$. We can write
$R=-E_1+\sum_{i=2}^tc_iE_i$.
But we have 
$$R=\gamma^{-1}(K_Y)+W-\tilde{p}^{-1}(p^{-1}(K_Y))-
\tilde{p}^{-1}(rF)=W-\tilde{p}^{-1}(rF).$$
Therefore we have $b_1=r-1$ and $c_i=b_i-ra_i$, for all $i\geq 2$.
It follows from the definition and from \cite{kollar},
Corollary 3.12, that the pair $(B,F)$ is canonical if and only if
$c_i\geq 0$ for all $i\geq 2$.

Since $B$ is locally a complete intersection, hence Gorenstein,
a result of Stevens \cite{stevens}
(see also\cite{kollar}, Theorem 7.9) says that $F$ is canonical
if and only if the pair $(B,F)$ is canonical near $F$. Since
$B\setminus F$ is smooth, this means precisely that $(B,F)$
is canonical.

On the other hand, since $F$ is locally a product of $X$
and an affine space, it follows that $X$ is canonical if and only if
$F$ is canonical, which completes the proof of the theorem.
\end{proof}

\begin{remark}\label{rational}
Since $X$ is locally a complete intersection, in particular Gorenstein,
it is a result of Elkik \cite{elkik} and Flenner \cite{flenner}
(see also \cite{komo},
Corollary~5.24) that $X$ has canonical singularities if and only if
it has rational singularities.
\end{remark}

We give also a necessary condition for $X$ to have
log canonical singularities and show that the sufficiency would follow
from the Inverse
of Adjunction Conjecture (\cite{kollar}, Conjecture 7.3).

\begin{theorem}\label{num_log_can}
With the above notation, if $b_i\geq ra_i-1$ for every $i$, then
$X$ has log canonical singularities. 
\end{theorem}

\begin{proof}
As in the proof of Theorem~\ref{num_can},
 we have $X$ log canonical if and only if
$F$ is log canonical and on the other hand $b_i\geq ra_i-1$ for all $i$
if and only if the pair $(B,F)$ is log canonical.

It follows from Proposition 7.3.2 in \cite{kollar} 
 that if $(B,F)$ is log canonical, then
$F$ is log canonical. 
\end{proof}

We conjecture that the converse is also true.

\begin{conjecture}\label{converse}
If $X$ has log canonical singularities, then 
$b_i\geq ra_i-1$, for all $i$.
\end{conjecture}

The argument in the proof of Theorem~\ref{num_log_can} shows that
Conjecture~\ref{converse} is implied by the conjecture below.
In fact, we will prove in the next section that these two conjectures
are equivalent.

\begin{conjecture}[Inverse of Adjunction]\label{inverse}
Let $X$ be a normal, l.c.i. variety and $D$ a normal Cartier divisor
on $X$. If $D$ is log canonical, then $(X,D)$ is log canonical
around $D$.
\end{conjecture}

\begin{remark}
In fact, the Inverse of Adjunction Conjecture is more general:
it deals with arbitrary normal varieties and with restriction of pairs
(see \cite{kollar}, Conjecture 7.3). It is known that it is implied
by the Log Minimal Model Program (see \cite{kollar_et_al},
Corollary 17.12). In particular, Conjecture~\ref{converse}
is true when $\dim\,Y=3$.
\end{remark}

\vfill\eject

\section{Irreducibility of Jet Schemes via Motivic Integration}

In this section we use motivic integration to give necessary
and sufficient conditions for a variety to have all the jet schemes
of the expected dimension and precisely one component of maximal dimension.
When the variety is locally a complete intersection, this gives via the
results in the previous two sections a proof of Theorem~\ref{main_conj}.

The construction of motivic integrals for smooth spaces
is due to Kontsevich \cite{kontsevich}, who used it to prove a conjecture
of Batyrev and Dais \cite{batdais} about the stringy Hodge numbers of
varieties with mild Gorenstein singularities
(see \cite{batyrev}). An other application
is the proof due to Batyrev \cite{batyrev2} of a conjecture 
of Reid on the McKay correspondence (see \cite{reid}).
The construction was generalized by 
 Denef and Loeser in \cite{denef1} and \cite{denef2}
to singular spaces (see the recent surveys
\cite{denef3} and \cite{looijenga} for other applications of this idea).
 We will need only
the Hodge realizations of motivic integrals on the space
of arcs of a smooth variety. We refer for definitions
and proofs to \cite{batyrev} (see also \cite{craw} for a nice
introduction).

{}From now on, $X$ will be a fixed 
variety over $k$, with ${\rm char}\,k=0$. Unless explicitely mentioned,
$X$ is not assumed to be locally complete intersection.
We fix an embedding 
$X\hookrightarrow Y$, where $Y$ is a smooth variety,
and an embedded resolution of singularities
$\gamma\,:\,\widetilde{Y}\longrightarrow Y$ for $(Y,X)$, 
as in the previous section.

More precisely,
we assume that $\gamma$ is a proper morphism which is an isomorphism
over $Y\setminus X$, and $\widetilde{Y}$ is smooth and $\gamma^{-1}(X)
=\sum_{i=1}^ta_iE_i$ is a divisor with normal crossings. Let 
$W=\sum_{i=1}^tb_iE_i$ be the discrepancy of $\gamma$. 
We set $N=\dim\,Y$ and $r={\rm codim}(X,Y)\geq 1$.
We can further assume that
$E_1$ is the only prime divisor
in $\gamma^{-1}(X)$ dominating $X$ and that $a_1=1$ and $b_1=r-1$.
With this notation, we prove the following results.

\begin{theorem}\label{num_expected_crit}
The following statements are equivalent:
\item{i)} For every $i\geq 1$, we have $b_i\geq ra_i-1$.
\item{ii)} $\dim\,X_m=(m+1)\dim\,X$, for every $m\geq 1$.
\item{iii)} There is $q\geq 1$, with $a_i\vert (q+1)$, for all $i$,
such that $\dim\,X_q=(q+1)\dim\,X$.
\end{theorem}

\begin{theorem}\label{num_irred_crit}
The following statements are equivalent:
\item{i)} For every $i\geq 2$, we have $b_i\geq ra_i$.
\item{ii)} For every $m\geq 1$, we have $\dim\,X_m=(m+1)\dim\,X$,
and $X_m$ has only one irreducible component of maximal dimension.
\item{iii)} There is $q\geq 1$, with $a_i\vert (q+1)$, for all $i$,
such that $\dim\,X_q=(q+1)\dim\,X$, and such that
$X_q$ has only one irreducible component of maximal dimension.
\end{theorem}

By combining Theorems~\ref{num_can} and \ref{num_irred_crit} and
Theorems~\ref{num_log_can} and \ref{num_expected_crit}, we obtain
the main results in our paper.

\begin{theorem}\label{main1}
If $X$ is an l.c.i. variety over $k$, and
${\rm char}\,k=0$, then $X_m$ is irreducible for every $m$
if and only if $X$ has canonical singularities.
\end{theorem}

\begin{proof}
Notice that by Proposition~\ref{normal} and by the definition
of canonical singularities, either condition implies that $X$ is normal,
so that Theorem~\ref{num_can} applies. Moreover, by Proposition~\ref
{irred_crit}, if $\dim\,X_m=(m+1)\dim\,X$, then $X_m$ is pure dimensional,
so that an application of Theorems~\ref{num_can}
and~\ref{num_irred_crit} completes the proof.
\end{proof}

Note that because of Remark~\ref{rational}, the above result is equivalent with
Theorem~\ref{main_conj}.

\begin{theorem}\label{main2}
If $X$ is a normal, l.c.i. variety over $k$ and
${\rm char}\,k=0$, and if $X_m$ is pure dimensional for every $m\geq 1$,
then $X$ has log canonical singularities. 
\end{theorem}

\begin{proof}
Again, Proposition~\ref{irred_crit} shows that $X_m$ is pure dimensional
if and only if $\dim\,X_m=(m+1)\dim\,X$ and we apply Theorems~\ref
{num_expected_crit} and~\ref{num_log_can}.
\end{proof}

\begin{conjecture}\label{converse2}
If $X$ is an l.c.i. variety over a field of
characteristic zero and $X$ has log canonical singularities, then 
$X_m$ is pure dimensional for every $m$.
\end{conjecture}

In fact, the above conjecture is equivalent with the conjectures we made
in the previous section.

\begin{proposition}\label{converse_and_inverse}
Conjectures~\ref{converse}, \ref{inverse} and \ref{converse2}
are equivalent.
\end{proposition}

\begin{proof}
Theorem~\ref{main1} implies that Conjectures~\ref{converse} and
\ref{converse2} are equivalent, and we have seen in the previous section
that Conjecture~\ref{inverse} implies Conjecture~\ref{converse}.
It is therefore enough to prove that if Conjecture~\ref{converse2}
is true for all normal, l.c.i. varieties, then so is Conjecture~\ref
{inverse}.
 
Using a trick due to Manivel (see \cite{kollar}, Lemma~7.1.3), the assertion
in Conjecture~\ref{inverse} can be reduced to the following: if $X$
is a normal, l.c.i.~variety and $D$ is a normal, Cartier divisor on $X$
which is log canonical, then $X$ is log canonical around $D$.

Applying Conjecture~\ref{converse2} to $D$, we get 
$\dim\,D_m=(m+1)(\dim\,X-1)$, for all $m$. We may assume that 
$D\hookrightarrow X$ is defined by one equation, so that $D_m
\hookrightarrow X_m$ is defined by $m+1$ equations. Therefore,
if $F$ is an irreducible component of $X_m$, with $\dim\,F>(m+1)\dim\,X$,
 then $F\cap D_m=\emptyset$. Since by Lemma~\ref
{irred_comp} we have $s_m(\overline{\pi_m(F)})\subseteq F$,
we deduce that $\overline{\pi_m(F)}\cap D=\emptyset$.
We conclude that there is an open neighbourhood $U$ of $D$ such that
$\dim\,U_m=(m+1)\dim\,U$. By Theorem~\ref{num_expected_crit},
it is possible to pick $m$ (depending on an embedded resolution of $X$)
such that if $\dim\,U_m=(m+1)\dim\,X$, then $U$ is log-canonical,
which finishes the proof of the proposition.
\end{proof}

By Lefschetz's principle, in order to prove Theorems~\ref{num_expected_crit}
and \ref{num_irred_crit}, we may assume $k=\CC$.  
Before we give the proof of the theorems, we review
the basics of motivic integration. For every $m\geq 0$,
we have natural projections
$$Y_{\infty}\overset{\eta_m}{\longrightarrow}
Y_m\overset{\pi_m}{\longrightarrow}Y,$$
where $\pi_m$ is an affine bundle with fiber $\AA^{mN}$.

The theory provides an
 algebra $\cM$ of subsets of $Y_{\infty}$, containing 
the algebra $\mathcal Cyl$ of cylinders of the form $\eta_m^{-1}(Z)$,
where $Z\subseteq Y_m$ is a constructible subset.
We consider the ring of Laurent power series in two variables
$u^{-1}$ and $v^{-1}$: $S=\ZZ[[u^{-1},v^{-1}]][u,v]$, with the
linear topology given by the descending sequence of subgroups
$\{\oplus_{i+j\geq l}\ZZ u^{-i}v^{-j}\}_l$. On $\cM$ there is a 
finitely additive measure $\mu$ with values in $S$, whose restriction to
$\mathcal Cyl$ is defined as follows.
If $C=\eta^{-1}_m(Z)$ is a cylinder, then
$$\mu(C)=E(Z;u,v)(uv)^{-(m+1)N}\in S,$$
where $E(Z;u,v)$ is the Hodge-Deligne polynomial of $Z$.
For a variety $Z$,
$$E(Z;u,v)=\sum_{1\leq p,q\leq\dim Z}\sum_{k\geq 0}
(-1)^kh^{p,q}(H^k_c(Z;\CC))u^pv^q,$$
where $\{h^{p,q}(H^k_c(Z;\CC))\}$ are the Hodge -Deligne numbers
of $Z$. What is important for us is that
$E(Z;u,v)$ is a polynomial of degree $2(\dim\,Z)$, and the
term of degree $2(\dim\,Z)$ is $c(uv)^{\dim\,Z}$, where
$c$ is the number of irreducible components of $Z$
of maximal dimension.

If $T\subseteq Y_{\infty}$ is a subset such that there is 
a sequence of cylinders $C_i$, with $T\subseteq C_i$
and $\mu(C_i)\longrightarrow 0$, then $T$ is in $\cM$ and
$\mu(T)=0$.

If $F\,:\,Y_{\infty}\longrightarrow\N\cup\{\infty\}$
is a function such that $F^{-1}(m)\in{\cM}$
for all $0\leq m\leq\infty$ and $\mu(F^{-1}(\infty))=0$,
and such that the series
$$\sum_{m\in\N}\mu(F^{-1}(m))(uv)^{-m}$$
is convergent in $S$, then $F$ is called integrable.
In this case, the sum of the above series is called
the motivic integral of $F$ and is denoted by
$\int_{Y_{\infty}}e^{-F}$.

In general, a subscheme $Z\subset Y$, with corresponding ideal
${\mathcal I}_Z$, defines a function $F_Z\,:\,Y_{\infty}
\longrightarrow\N\cup\{\infty\}$, as follows. If $w\in 
Y_{\infty}$ is an arc over $y\in Y$, then $w$ can be
identified with a ring homomorphism $w\,:\,{\mathcal O}_{Y,y}
\longrightarrow\CC[[t]]$ and $$F_Z(w)={\rm ord}
(w({\mathcal I}_{Z,y})).$$
In \cite{batyrev} and \cite{craw}, the authors consider this function when
$Z$ is a divisor in $Y$. It follows from the definition that
$F_Z^{-1}(\infty)=Z_{\infty}$ and that
 $$F_Z^{-1}(m)=\eta_{m-1}^{-1}(Z_{m-1})\setminus\eta_m^{-1}(Z_m),$$
for any integer $m\geq 0$ (we make the convention that 
$Z_{-1}=Y$ and $\eta_{-1}=\eta_0$). In particular, $F_Z^{-1}(m)$
is a cylinder.

\begin{lemma}\label{upper_bound}
If $D\subset Y$ is a divisor and $y\in Y$ is a point such that
${\rm mult}_yD=a$, then $\dim\,{(\pi^D_m)}^{-1}(y)
\leq Nm-[m/a]$, where $[x]$ denotes the integral part of $x$.
In particular, for every proper subscheme $Z\subset Y$,
we have $\dim\,Z_m-N(m+1)\longrightarrow -\infty$.
\end{lemma}

\begin{proof}
The second assertion follows from the first one, since working locally
we may assume that $Z\subset D$, for some divisor $D$ and if
$a=\max_{y\in Y}\{{\rm mult}_yD\}$, then
$\dim\,Z_m-N(m+1)\leq -[m/a]$.

To prove the first statement, note that the case $a=0$ is trivial, and
 therefore we may assume that $y\in D$. It is enough to show that for every 
$p\geq 1$ we have
\begin{equation}\label{ineq}
\dim\,(\pi_{pa}^D)^{-1}(y)\leq paN-p.
\end{equation}

If we pick a regular system of parameters $x_1\ldots,x_N$ in ${\mathcal O}
_{Y,y}$, we get an \'{e}tale ring homomorphism 
$w\,:\,{\mathcal O}_{{\bf A}^N,0}\longrightarrow
{\mathcal O}_{Y,y}$. Let $f\in {\mathcal O}_{Y,y}$ be an equation for
$D$ at $y$. In general $f\not\in{\rm Im}\,(w)$. However, if 
$g\in {\mathcal O}_{Y,y}$ is such that $g-f\in \underline{m}_{Y,y}^{pa+1}$,
then $(\pi_{pa}^D)^{-1}(y)=(\pi_{pa}^{V(g)})^{-1}(y)$. Since
$w$ induces an isomorphism of the associated graded rings, we can
find $g$ as before such that $g\in {\rm Im}\,(w)$. Therefore, in order
to prove (\ref{ineq}), we may assume that $f\in {\rm Im}\,(w)$.
Since $w$ is \'{e}tale, by replacing ${\bf A}^N$ and $Y$ with
suitable open neighbourhoods of $0$ and $y$, respectively, 
 we can apply Lemma~\ref
{etal_com} to reduce to the case when $Y={\bf A}^N$, $y=0$ and
$D$ is defined by a polynomial $f\in k[X]=k[X_1,\ldots, X_N]$.

 We  use the equations for $D_{pa}$ described in the first
section. With the notation
 $f^{(j)}_0=f^{(j)}(0,X',\ldots,X^{(j)})$
for every $j$, we have 
$(\pi_{pa}^D)^{-1}(0)\subseteq 
Z_p$, where $Z_p$ is defined in
$\Spec\,k[X',\ldots,X^{(pa)}]$ by the polynomials 
$f^{(ja)}_0$, with $1\leq j\leq p$.

We prove that $\dim\,Z_p\leq paN-p$ by computing the dimension of
 a deformation of this set. Note that if we put $\deg(X_i^{(j)})=j$ for
every $i$ and $j$, then each polynomial $f^{(j)}_0\in k[X',\ldots,X^{(pa)}]$
is homogeneous of degree $j$. 

Consider the family ${\mathcal Z}_p$ over $\Spec\,k[t]$ defined
in ${\bf A}^{paN}\times \Spec\,k[t]$ by the
polynomials $(1/t^a)f^{(j)}_0(tX',\ldots,tX^{(j)})$, with $1\leq j\leq
p$. The fiber of ${\mathcal Z}_p$ over every $t_0\neq 0$ is isomorphic to 
$Z_p$, while the fiber over $0$ is the corresponding scheme 
obtained by replacing $f$ with
its homogeneous component of degree $a$. Since 
all the rings are graded (for the grading we defined above), the 
semicontinuity theorem for the dimension of the fibers of a morphism
(see \cite{eisenbud}, Theorem 14.8)
shows that in order to prove that $\dim\,Z_p\leq paN-p$, we may assume
that $f$ is homogeneous of degree $a$.

Consider now on $k[X,X',\ldots,X^{(pa)}]$ the reverse lexicographic order
where we order the variables such that $X_i^{(j)}<X_{i'}^{(j')}$
if $j>j'$ or if $j=j'$ and $i>i'$ (see, for example, 
\cite{eisenbud} Chapter 15). Let $m(X)={\rm in}(f)$ be the initial term of
$f$ is this order. It is then easy to see that ${\rm in}(f^{(ja)}_0)
=m(X^{(j)})$, for $1\leq j\leq p$. Therefore the initial ideal of
the ideal defining $Z_p$ has dimension $paN-p$.
Since the dimension of an ideal is equal with the
dimension of its initial ideal, we deduce that 
$\dim\,Z_p=paN-p$, which concludes the proof of the lemma.
\end{proof}

\begin{corollary}\label{measurable}
For every proper subscheme $Z\subset Y$, and every $0\leq m\leq\infty$,
$F_Z^{-1}(m)\in\cM$ and $\mu(F_Z^{-1}(\infty))=0$.
\end{corollary}

\begin{proof}
We have already seen that $F_Z^{-1}(m)$ is a cylinder for every integer
$m\geq 0$.
Moreover,  $F_Z^{-1}(\infty)=Z_{\infty}\subseteq\eta_m^{-1}(Z_m)$ and
Lemma~\ref{upper_bound} gives
$$\mu(\eta_m^{-1}(Z_m))=E(Z_m;u,v)(uv)^{-N(m+1)}\longrightarrow 0.$$
\end{proof}

In order to prove Theorems~\ref{num_expected_crit} and \ref{num_irred_crit},
we will choose a suitable function $f\,:\,\N\longrightarrow\N$
(which we extend by $f(\infty)=\infty$) and we will integrate
$F=f\circ F_X$ on $Y_{\infty}$.
 The change of variable
formula (see \cite{batyrev}, Theorem 6.27 or \cite{denef1}, Lemma 3.3)
gives
$$\int_{Y_{\infty}}e^{-F}=\int_{\widetilde{Y}_{\infty}}
e^{-(F\circ \gamma_{\infty}+F_W)},$$
in the sense that one integral exists if and only if the other one does,
and in this case they are equal. The point is that 
$F\circ \gamma_{\infty}=f\circ F_{\gamma^{-1}(X)}$,
 and since $\gamma^{-1}(X)\cup W$
has normal crossings, the right-hand side integral can be explicitely
computed, while for a suitable choice of $f$, the left-hand side
contains the information we need about the dimension of $X_m$ and about the
number of its irreducible components of maximal dimension.

\begin{proof}[Proof of Theorems~\ref{num_expected_crit} and
~\ref{num_irred_crit}]

We fix a function $f\,:\,\N\longrightarrow\N$, such that for every
$m\geq 0$,
$${(\star)}\, f(m+1)>f(m)+\dim\,X_m+C(m+1),$$
where $C\in\N$ is a constant with
$C>|N-(b_j+1)/a_j|$, for all $j$. We extend it by defining 
$f(\infty)=\infty$. For the proof of the implication $iii)\Rightarrow i)$
we will put later an extra condition.

It follows from Corollary~\ref{measurable} that if $F=f\circ F_X$, then
$F^{-1}(m)\in{\mathcal M}$ for $0\leq m\leq\infty$, and 
$\mu(F^{-1}(\infty))=0$.  Computing the integral of $F$ from the
definition, we get $I=\int_{Y_{\infty}}e^{-F}=S_1-S_2$, where
$$S_1=\sum_{m\geq 0}E(X_{m-1};u,v)(uv)^{-mN-f(m)},$$
$$S_2=\sum_{m\geq 0}E(X_m;u,v)(uv)^{-(m+1)N-f(m)}.$$

Every monomial which appears in the $m^{\rm th}$ term of $S_1$,
has degree bounded above by $2P_1(m)$ and below by $2P_2(m)$,
where 
$$P_1(m)=\dim\,X_{m-1}-mN-f(m),$$ 
$$P_2(m)=-mN-f(m),$$ for all $m\geq 0$
(recall the convention that $X_{-1}=Y$). 
Moreover, we always have precisely one monomial of degree $2P_1(m)$,
namely $(uv)^{P_1(m)}$, whose coefficient is $c_m$, the number of irreducible
components of maximal dimension of $X_{m-1}$.

Similarly, every monomial which appears in the $m^{\rm th}$ term of
$S_2$ has degree bounded above by $2Q_1(m)$ and below by $2Q_2(m)$,
where 
$$Q_1(m)=\dim\,X_m-(m+1)N-f(m),$$
$$Q_2(m)=-(m+1)N-f(m),$$
for all $m\geq 0$. 
We always have exactly one
 monomial of degree $2Q_1(m)$, namely $(uv)^{Q_1(m)}$,
whose coefficient is $c_{m+1}$.

A first consequence of this and Lemma~\ref{upper_bound} is that $F$
is, indeed, integrable. Using condition ${(\star)}$, it is an easy
 computation to show that we have $P_1(m+1)<\min\{P_2(m), Q_2(m)\}$,
for every $m\geq 0$.

Moreover, Lemma~\ref{fibers} gives $\dim\,X_m\leq\dim\,X_{m-1}+N$, for every
$m\geq 1$, and Lemma~\ref{upper_bound} implies that the inequality is strict 
for infinitely many $m$. We deduce that $Q_1(m)\leq P_1(m)$, for every
$m\geq 0$ and equality holds if and only if $m\geq 1$ and
$\dim\,X_m=\dim\,X_{m-1}+N$ (therefore, the inequality is strict for
infinitely many $m$). 

We conclude from the above inequalities first that in $S_1$,
the term 
$(uv)^{P_1(m)}$ appears precisely once for every $m\geq 0$, 
and has coefficient
$c_m$. Similarly, in $S_2$, the term $(uv)^{P_1(m)}$ appears at most once. 
It appears if and only if $m\geq 1$ and $\dim\,X_m=\dim\,X_{m-1}+N$, and 
in this case it has coefficient $c_{m+1}$.

We use the change of variable formula to compute the integral of $F$,
as
 $$I=\int_{Y_{\infty}}e^{-F}=
\int_{\widetilde{Y}_{\infty}}e^{-(f\circ F_{\gamma^{-1}(X)}+F_W)}.$$
In this form, $I$ can be explicitely computed, since $\gamma^{-1}(X)\cup W$ has
normal crossings. For every subset $J\subseteq\{1,\ldots,t\}$,
let $E_J^{\circ}=\cap_{i\in J}E_i\setminus\cup_{i\not\in J}E_i$.
With this notation,
we have 

$$\int_{\widetilde{Y}_{\infty}}e^{-(f\circ F_{\gamma^{-1}(X)}+F_W)}=
\sum_{J\subseteq\{1,\ldots\,t\}}S_J,\,{\rm where}$$

$$S_J=\sum_{\alpha_i\geq 1,i\in J}E(E_J^{\circ};u,v)
(uv-1)^{|J|}\cdot (uv)^{-N-\sum_{i\in J}\alpha_i(b_i+1)-f(\sum_{i\in J}
a_i\alpha_i)}.$$
We just sketch the proof of this formula, as it is similar
to that of Theorem~6.28 in \cite{batyrev}, or that
of Theorem~1.16 in \cite{craw}.

Since $E$ is additive, we may work locally on $\widetilde{Y}$.
In order to compute the part $S_J$ in $I$ which corresponds
to arcs over $E_J^{\circ}$, we may assume that there is
a regular system of parameters $y_1,\ldots,y_N$ on $Y$, such
that $E_i$ is defined by $y_i$, for all $i\in J$. We have
$$S_J=\sum_{\alpha_i\geq 1,i\in J}\mu(\{u\in\eta_0^{-1}(E_J^{\circ})
\colon F_{E_i}(u)=\alpha_i,i\in J\})(uv)^{-f(\sum_{i\in J}a_i\alpha_i)-
\sum_{i\in J}b_i\alpha_i}.$$
For every $(\alpha_i)_{i\in J}$, if $p\geq {\rm max}_{i\in J}\{\alpha_i\}$,
then
$$\{u\in\eta_0^{-1}(E_J^{\circ})\colon F_{E_i}(u)=\alpha_i, i\in J\}=
\eta_p^{-1}(C^p_{\alpha_i,i\in J}),$$
where $C^p_{\alpha_i,i\in J}\hookrightarrow Y_p$ is locally trivial
over $E_J^{\circ}$, with fiber ${(k^*)}^{|J|}\times\AA^{Np-\sum_{i\in J}
\alpha_i}$, and our formula for $S_J$ follows.

Every monomial in the term of $S_J$ corresponding to $(\alpha_i)_{i\in J}$,
has degree bounded above by $2R_1(\alpha_i;i\in J)$, and below
by $2R_2(\alpha_i;i\in J)$, where 
$$R_1(\alpha_i;i\in J)=-\sum_{i\in J}\alpha_i(b_i+1)-f(\sum_{i\in J}
a_i\alpha_i)$$
and $R_2(\alpha_i;i\in J)=R_1(\alpha_i;i\in J)-N$. Note that if
$J=\emptyset$, then $R_1(\emptyset)=-f(0)$.

We introduce one more piece of notation: for $m\geq 0$, let
$\tau(m)=\dim\,X_m-(m+1)\dim\,X$. For every $m\geq 0$, we have $\tau_m\geq 0$.
We see that for $J\neq\emptyset$, we have
$$R_1(\alpha_i,i\in J)=P_1(\sum_{i\in J}a_i\alpha_i)
-\tau(\sum_{i\in J}a_i\alpha_i-1)-\sum_{i\in J}\alpha_i(b_i+1-ra_i).$$

Moreover, property ${(\star)}$ implies that if $J\neq\emptyset$, then
\begin{equation}\label{eq1}
P_1(\sum_{i\in J}a_i\alpha_i+1)<
R_2(\alpha_i;i\in J)<R_1(\alpha_i;i\in J)<
\end{equation}
$$\min\{P_2(\sum_{i\in J}
a_i\alpha_i-1), Q_2(\sum_{i\in J}a_i\alpha_i-1)\}$$
and that $P_1(1)<R_2(\emptyset)$. In particular,
this implies that the only monomial of the form $(uv)^{P_1(m)}$
which can appear in the term corresponding to
$J$ and  $(\alpha_i)_{i\in J}$ is for
$m=\sum_{i\in J}a_i\alpha_i$.

To prove the implication $i)\Rightarrow ii)$ in Theorem~\ref{main1},
suppose that $b_i\geq ra_i-1$, for all $i$ and assume that for some
$m\geq 1$, we have $\tau(m)>0$. The above inequalities show that
$(uv)^{P_1(m+1)}$ does not appear in the sum $S_J$, for
every $J$.

As we have seen, this imples that $\dim\,X_{m+1}=\dim\,X_m+N$. 
In particular, we have $\tau(m+1)>0$. Continuing in this way, we get
$\dim\,X_{p+1}=\dim\,X_p+N$, for every $p\geq m$, a contradiction with
Lemma~\ref{upper_bound}. Therefore we must have $\tau(m)=0$, for all
$m$, so $\dim\,X_m=(m+1)\dim\,X$.

Suppose next that $b_i\geq ra_i$, for all $i\geq 2$. The above argument 
shows that $\tau(m)=0$, for every $m\geq 0$. In particular, the coefficient
of $(uv)^{P_1(m+1)}$ in $I$ is $c_{m+1}$, for every $m\geq 0$.  
{}From the above inequalities, we see that for every $m\geq 0$,
the term 
$(uv)^{P_1(m+1)}$ appears in $S_J$ if and only if
 $J=\{1\}$ and in this
case it has coefficient $1$, since $E_{\{1\}}^{\circ}$
 is irreducible. Therefore
$c_{m+1}=1$, for every $m\geq 0$. 
This proves the implication $i)\Rightarrow ii)$ in Theorems
~\ref{num_expected_crit} and~\ref{num_irred_crit}.

 We next turn to 
the implication $iii)\Rightarrow i)$. Suppose that 
for some $q\geq 1$, with $a_i\vert (q+1)$, for all $i$, we have $\tau(q)=0$,
 and that for some $j\leq t$, $b_j<ra_j-1$.

We pick the function $f$ such that in addition to $(\star)$ it satisfies 
the following requirement. For every $p$, consider the set 
${\mathcal J}_p$ of all the pairs $(J, (\alpha_i)_{i\in J})$, such that
$\sum_{i\in J}a_i\alpha_i=p$. This is clearly a finite set. We require
that for every $(J, (\alpha_i)_{i\in J})\in {\mathcal J}_{q+1}$
 and every $(J', (\alpha'_i)_{i\in J'})\in
{\mathcal J}_p$, for some $p$ with $p\leq q$, we have
\begin{equation}\label{requirement}
f(q+1)>f(p)-\sum_{i\in J}\alpha_i(b_i+1)+\sum_{i\in J'}\alpha'_i(b_i+1)+N,
\end{equation}
or equivalently, we have $R_1(\alpha_i;i\in J)<R_2(\alpha'_i;i\in J')$.

On the other hand, if $(J',(\alpha'_i)_{i\in J'})\in {\mathcal J}_p$, for some
$p\geq q+2$, then by $(\ref{eq1})$ we have
\begin{equation}\label{eq2}
 P_2(q+1) > R_1(\alpha_i';i\in J').
 \end{equation}

Note that the top degree monomials 
which appear in   different terms of the sums
$S_J$ (for possibly different $J$)
 don't cancel each other, because they have positive
coefficients. Let $d$ be the highest degree of a monomial which appears
in a term corresponding to some $(J, (\alpha_i)_{i\in J})\in
 {\mathcal J}_{q+1}$. By the previous remark, the corresponding monomial
does not cancel with a monomial in a term corresponding to
$(J',(\alpha'_i)_{i\in J'})\in {\mathcal J}_{q+1}$.

Since $(\{j\}, (q+1)/a_j)\in {\mathcal J}_{q+1}$, our hypothesis implies that
$2P_1(q+1)<d$, while from $(\ref{eq1})$ we deduce $d<2\min\{P_2(q), Q_2(q)\}$.
Moreover, we deduce from
$(\ref{requirement})$ and $(\ref{eq2})$ that the monomial of degree $d$
 does not cancel
with monomials in terms
 corresponding to $(J', (\alpha'_i)_{i\in J'})\in {\mathcal J}_p$,
if $p\neq q+1$. This shows that in $I$ there is
indeed a monomial of degree $d$,
where $2P_1(q+1)<d<2\min\{P_2(q), Q_2(q)\}$, a contradiction.

Suppose next that for some $q\geq 1$,
with $a_i\vert (q+1)$, for all $i$,
 we have $\tau(q)=0$ and $c_{q+1}=1$,
and that for some $j\geq 2$, $b_j<ra_j$. The above argument shows that
$b_i\geq ra_i-1$, for every $i$. Note that since 
$a_j\vert (q+1)$, we have in the expression of
 $I$ the monomial $(uv)^{P_1(q+1)}$, with coefficient
at least $2$, once from $S_{\{1\}}$ (with $\alpha_1=
q+1$), and once from $S_{\{j\}}$ (with
$\alpha_j=(q+1)/a_j$). This gives a contradiction.
\end{proof}

\section{Examples and Open Problems}

In this section, unless explicitely mentioned otherwise, $k$
has arbitrary characteristic.

We consider first the case of curves and show that
none of the higher jet schemes of a singular curve
can be irreducible.  Since we do not assume
$X$ to be locally complete intersection, we need first a lemma showing
that if the tangent space at a point
to a scheme is too big, then so are all the 
fibers over that point of the higher order jet schemes.

\begin{lemma}\label{dim_TX}
If $X$ is a scheme and $x\in X$, then
$\dim\,\pi_{2m}^{-1}(x)\geq \dim_xX+m \dim\,T_xX$ and
$\dim\,\pi_{2m+1}^{-1}(x)\geq (m+1)\dim T_xX$, for all $m\geq 1$.
\end{lemma}

\begin{proof}
We prove the first assertion. Let $f\,:\,\pi_{2m-1}^{-1}(x)
\longrightarrow \pi_{m-1}^{-1}(x)$ be the canonical projection.
Using Lemma~\ref{fibers}, it is enough to show that
there is an isomorphism $T_xX^{\oplus m}\simeq f^{-1}(0)$,
which maps $C_xX\times T_xX^{\oplus(m-1)}
\subseteq T_xX^{\oplus m}$
into $f^{-1}(0)\cap {\rm Im}(\phi_{2m})$. Here $0$ denotes the image of $x$
by the canonical section of $\pi_{m-1}$ and $C_xX$ is the tangent cone to
$X$ at $x$.

We give the isomorphism at the level of $A$-valued points.
An $A$-point of $f^{-1}(0)$ is given by an algebra homomorphism
$\theta\,:\,\cO_{X,x}\longrightarrow A[t]/(t^{2m})$ of the form
$\theta(y)=\theta_0(y)+\sum_{i=m}^{2m-1}\theta_i(y)t^i$, where
$\theta_0$ corresponds to $x$. The condition that $\theta$ is an
algebra homomorphism is equivalent with saying that for
every $i$, with $m\leq i\leq 2m-1$, the morphisms $\theta'_i$
mapping $y$ to $\theta_0(y)+\theta_i(y)t$ are $A$-valued points of
$T_xX$.

If $\theta$ is a $k$-valued point of $f^{-1}(0)$, then $\theta\in {\rm Im}
(\phi_{2m})$ if and only if the morphism
$S^2(m_x/m_x^2)\longrightarrow k$, given by $y_1\cdot y_2\longrightarrow
\theta_m(y_1)\theta_m(y_2)$ factors through $m_x^2/m_x^3$. This condition
is satisfied in particular if $\theta'_m$ is a $k$-valued point of $C_xX$.

The proof of the second assertion is similar, giving for
the projection $g\,:\,\pi_{2m+1}^{-1}(x)\longrightarrow\pi_m^{-1}(x)$,
an isomorphism $g^{-1}(0)\simeq T_xX^{\oplus(m+1)}$.
\end{proof}

\begin{corollary}
If $X$ is an integral curve,
 then for any $m\geq 1$, $X_m$ is irreducible if and only if
$X$ is nonsingular.
\end{corollary}

\begin{proof}
Indeed, if $x\in X$ is a singular point, then $\dim\,T_xX\geq 2$,
and by Lemma~\ref{dim_TX}, it follows that for every $m\geq 1$, 
$\dim\pi_m^{-1}(x)\geq m+1$, and therefore $\pi_m^{-1}(x)$
gives an irreducible component of $X_m$.
\end{proof}

\begin{remark}
If ${\rm char}\,k\ne 2$, it is possible to show that
an integral curve $X$ has all the jet schemes pure dimensional if and only if
it has at worst nodes as singularities.
\end{remark}

We consider next the case of a surface.  In this case, if
all the jet schemes of $X$ are irreducible, 
we show that $X$ has to be locally
complete intersection, so we can apply our previous theory.

\begin{theorem}\label{surfaces}
Assume that ${\rm char}\,k=0$ and let $X$ be a surface. Then the 
following are equivalent:
\item{i)} $X$ has at worst rational double points as singularities.
\item{ii)} $X_m$ is irreducible for every $m\geq 1$.
\item{iii)} $X_m$ is irreducible for $m\gg 0$.
\end{theorem}

\begin{proof}
One of the characterizations of rational double points is that they are
locally complete intersection rational singularities
(see \cite{durfee}). Therefore 
``$i)\Rightarrow ii)$'' follows from Theorem~\ref{main1}
and the fact that Gorenstein singularities
are rational if and only if they
are canonical (see Remark~\ref{rational}).
In order to prove $iii)\Rightarrow i)$, it is enough to show
that if $X_m$ is 
irreducible for some $m$, then $X$ is locally complete intersection,
and then use Proposition~\ref{one_implies_lower} and Theorem~\ref{main1}.
But for every $x\in X$, Lemma~\ref{dim_TX} implies that if
$\dim T_xX\geq 4$, then $\dim\pi_m^{-1}(x)\geq 2m+2$, and therefore
$X_m$ is not irreducible. We conclude that $\dim T_xX=3$ at every
singular point $x\in X$, in particular that $X$ is locally complete
intersection. 
\end{proof}

\begin{remark}
It is possible to prove Theorem~\ref{surfaces} directly, by showing that
condition $ii)$ is equivalent to $X$ having the singular points of one
of the types in the classification of the rational double points. In fact,
that proof shows more, namely that the above conditions are equivalent
with having only $X_5$ irreducible.
\end{remark}

\begin{example}
Theorem~\ref{main1} is not true if we replace the condition
of $X$ being locally a complete intersection variety, with being Gorenstein.
More precisely, it is possible to have a variety with Gorenstein
 canonical singularities, but such that all its jet schemes are not even
pure dimensional.

In fact, it is possible to take $X$ to be a toric variety of dimension $3$
(for definitions, basic facts and notations for toric varieties, see 
\cite{fulton}). If $V$ is a $\QQ$-vector space with basis 
$e_i$, $1\leq i\leq 3$, let $N\subset V$ be the lattice spanned by
$\{e_1,e_2,e_3,1/3(e_1+e_2+e_3)\}$. If $\sigma$ is the cone in $V$
spanned by $\{e_i\}_i$, let $X=U_{\sigma}$ be the associated toric variety.
Then $X$ is Gorenstein, since $\sum_ie_i^*\in N^{\vee}$
 (see \cite{fulton}, Sections 3.4 and 4.4).
By \cite{fulton}, Section 3.5, it has rational,
hence canonical singularities (see Remark~\ref{rational}). Moreover, 
$\{\sum_ia_ie_i^*\,\mid\,a_i\in\NN,\sum_ia_i=3\}$ induce linearly independent
elements in $T_xX$, where $x\in X$ is the fixed point under the
torus action.
Therefore, $\dim T_xX\geq 10$, and by Lemma~\ref{dim_TX} we get
for every $m\geq 1$, $\dim\pi_m^{-1}(x)\geq 5m+3$, and since $\dim\,X=3$, 
$X_m$ is not pure dimensional.
\end{example}
 
\begin{example}
On the other hand, the condition that $X$ is locally a complete intersection
is not necessary in order to have $X_m$ irreducible for every $m$.

Let $X\subset\AA^{2n}$ be the cone over the Segre embedding 
$\PP^1\times\PP^{n-1}\hookrightarrow\PP^{2n-1}$. It is defined by the
ideal generated by the $2\times 2$ minors of the generic matrix:

$$\begin{pmatrix}
U_1&U_2&\ldots&U_n\\
V_1&V_2&\ldots&V_n\\ 
\end{pmatrix}$$
in the ring $S=k[U_1\ldots,U_n,V_1,\ldots,V_n]$. Notice that if $n\geq 3$,
then $X$ is not a complete intersection, but $X_m$
is irreducible for every $m\geq 1$.

Indeed, since $X$ is defined by degree two homogeneous polynomials, it is easy
to see that $\pi_m^{-1}(0)\cong X_{m-2}\times \AA^{2n}$, for all
$m\geq 1$ (we take $X_{-1}$ to be a point). By induction on $m$,
we get $\pi_m^{-1}(0)$ irreducible, and because $X_{\rm sing}=\{0\}$,
it is enough to find a nonempty subset $U\subset\pi_m^{-1}(0)$ such that
$U\subset\overline{\pi_m^{-1}(X_{\rm reg})}$. We consider the open subset
of $\pi_m^{-1}(0)$ given by matrices
\begin{equation*}
A=
\begin{pmatrix}
a_1&a_2&\ldots&a_n\\
b_1&b_2&\ldots&b_n\\
\end{pmatrix}
\end{equation*}
with $a_i$, $b_i\in k[t]/(t^{m+1})$ of order one, such that 
$a_ib_j=a_jb_i$, for all $i$ and $j$. There is a unique 
$p=p_0+p_1t+\ldots +p_{m-2}t^{m-2}$, $p_0\neq 0$, such that
for all $i$, $b_i=pa_i+c_it^m$, for some $c_i\in k$.

By restricting to a smaller open subset, we may assume that $c_i\neq 0$,
for all $i$. In this case $A\in\overline{\{A_s\mid s\in k^*\}}$, where
$A_s\in\pi_m^{-1}(X_{\rm reg})$ is given by
\begin{equation*}
A_s=
\begin{pmatrix}
a_1+c_1s&\ldots&a_n+c_ns\\
p_s(a_1+c_1s)&\ldots&p_s(a_n+c_ns)\\
\end{pmatrix}
\end{equation*}
where $p_s=p+1/s\cdot t^m$.
\end{example}

\begin{question}
What is the analogue of Theorem~\ref{main1} in positive characteristic ?
Is it true that if $X$ is an l.c.i. variety, then $X_m$ is irreducible for all
 $m$ if and only if $X$ has pseudorational singularities ? The notion
of pseudorational singularities, introduced by Lipman and Teissier in
\cite{lipman}, replaces the notion of rational singularities when 
a good desingularization theory and Grauert-Riemenschneider theorem
are not known. When these results are known (for example, in characteristic
zero or for surfaces), the two notions coincide.

A different analogue of rational singularities in positive characteristic,
coming from tight closure theory,
is that of F-rational singularities (see \cite{smith1}, for definition
and relations with the birational geometry). A result of Smith 
(see \cite{smith2}, Theorem 3.1) says that F-rational singularities
are pseudorational. In general, having F-rational singularities is not a
 necessary condition for having irreducible jet schemes. For example, it
follows from Proposition~\ref{quasihom} below, that $V(X^2+Y^3+Z^5)$
has irreducible jet schemes in any characteristic,
 while it is known that it is not
F-rational if ${\rm char}\,k\in\{2,3,5\}$ (see \cite{bruns}, 
Example 10.3.12).
\end{question}

\begin{proposition}\label{quasihom}
For $n\geq 3$, let $F=X_1^{d_1}+\ldots+X_n^{d_n}$ and $Z=V(F)
\subset\AA^n$, such that there is at most one $i$, with 
${\rm char}\,k|d_i$. We have $Z_m$ irreducible for all $m\geq 1$
if and only if $\sum_{i=1}^n1/d_i>1$.
\end{proposition}

\begin{proof}
The hypothesis says that $Z\setminus\{0\}$ is smooth; in particular,
$Z$ is integral. By Proposition~\ref{irred_crit}, $Z_m$ is irreducible 
for all $m\geq 1$ if and only if $\dim\,\pi_m^{-1}(0)
<(m+1)(n-1)$ for all $m$, and by Proposition~\ref{one_implies_lower},
it is enough to check this for infinitely many $m$.

For every $m\geq 1$ and integers $a_1,\ldots,a_n$, with
$1\leq a_i\leq m+1$ for all $i$, let $V_{a_1,\ldots,a_n}$
be the locally closed subset of $\pi_m^{-1}(0)$ consisting of
ring homomorphisms $\phi\,:\,k[X_1,\ldots,X_n]/(F)\longrightarrow
k[t]/(t^{m+1})$ with ${\rm ord}\,(\phi(X_i))=a_i$ (we make the
convention ${\rm ord}(0)=m+1$). We obviously have
$\pi_m^{-1}(0)=\bigcup_{a_1,\ldots,a_n}V_{a_1,\ldots,a_n}$.

Set $d=\prod_{i=1}^n
d_i$ and assume that $Z_{d-1}$ is irreducible.
If $a_i=d/d_i$, for every $i$, then $\phi(X_i)$ can be
chosen arbitrarily with order $a_i$, so that $\dim\,V_{a_1,\ldots,a_n}
=d(n-\sum_i1/d_i)$. Since $\dim\,V_{a_1,\ldots,a_n}<
d(n-1)$, we get $\sum_i1/d_i>1$. 

Next suppose that $\sum_i1/d_i>1$ and that  $d|m+1$.
We will show that
$Z_m$ is irreducible. Consider first the case when $d_ia_i\geq m+1$,
for all $i$. As above, in this case $\phi(X_i)$ can be chosen
arbitrarily with order $a_i$. Therefore we get
$$\dim\,V_{a_1,\ldots,a_n}=\sum_i(m-a_i+1)\leq n(m+1)-(m+1)\sum_i1/d_i
<(m+1)(n-1).$$

On the other hand, if
$r={\rm inf}_i\{a_id_i\}\leq m$, and
if $V_{a_1,\ldots,a_n}\neq\emptyset$, then by hypothesis there is $i$,
such that ${\rm char}\,k\not\ \mid d_i$, and $a_id_i=r$.
For every choice of $\phi(X_j)$ with order $a_j$, $j\neq i$, 
we have $\phi(X_i)=t^{a_i}u$, where the image of 
$u$ in $k[t]/(t^{m+1-r})$ can take finitely many values. Moreover,
$\phi(X_i)$ is uniquely determined by the image of $u$ in
$k[t]/(t^{m+1-a_i})$. Therefore we have
$$\dim\,V_{a_1,\ldots,a_n}\leq$$
$$ r-a_i+\sum_{j\neq i}
(m+1-a_j)\leq (n-1)(m+1)+r(1-\sum_{j=1}^n1/d_j)<(n-1)(m+1).$$
It follows that $\dim\,\pi_m^{-1}(0)<(n-1)(m+1)$, so $Z_m$ is irreducible.
Since this is true for all $m$ with $d\vert (m+1)$, it follows by
Proposition~\ref{one_implies_lower} that $Z_m$ is irreducible for all $m$.
\end{proof}

\begin{question}\label{reduced2}
We have studied when an l.c.i. variety has
irreducible jet schemes. Proposition~\ref{reduced} shows that if
$X_m$ is irreducible, then $X_m$ is reduced. Is the converse true ?
 We will see below that this the case for $m=1$.
\end{question}

\begin{question}
Is it true that if $X$ is an l.c.i. variety, then
for every $m\geq 1$, $(X_m)_{\rm reg}=\pi_m^{-1}(X_{\rm reg})$ ?
A positive answer to this question would give a positive answer
to Question~\ref{reduced2}. Indeed, if $X_m$ is reduced, but not
irreducible, then any irreducible component contained in
$\pi_m^{-1}(X_{\rm sing})$ has to be generically smooth.

We can give a positive answer when $m=1$. This is the content of the
following proposition.
\end{question}

\begin{proposition}\label{case1}
If $X$ is an l.c.i. variety, then $(X_1)_{\rm reg}
=\pi_1^{-1}(X_{\rm reg})$. 
\end{proposition}

\begin{proof}
We show first that if $u\in\pi_1^{-1}(x)$, then
$$(\dagger)\,\,\dim\,T_uX_1\geq\dim\,X+\dim\,T_xX.$$
To see this, we may work locally and assume that $X\subset\AA^n$ is defined by 
$\underline{f}=(f_1,\ldots,f_r)$, where $r=n-\dim\,X$.
$X_1$ is defined by $(\underline{f},\underline{f}')$, and
if $u=(x,x')\in X_1\subset\AA^{2n}$, then the Jacobian
$J_{(\underline{f},\underline{f}')}(u)$ is
$$\begin{pmatrix}
J_{\underline{f}}(x)&0\\
A&J_{\underline{f}}(x)\\
\end{pmatrix}$$
for some $r\times n$ matrix $A$. Therefore ${\rm rk}\,J_{(\underline{f},
\underline{f}')}(u)\leq r+{\rm rk}\,J_{\underline{f}}(x)$
and $(\dagger)$ follows.

Suppose now that $u\in (X_1)_{\rm reg}$.
Consider an open connected neighbourhood $U\subset (X_1)_{\rm reg}$ of $u$.
The inequality $(\dagger)$ for an arbitrary point $u'\in U$ gives
$$\dim\,U\geq\dim\,X+\dim\,\pi_1^{-1}(\pi_1(u')).$$ 

This implies that $\overline{\pi_1(U)}=X$, and therefore $U\cap
\pi_1^{-1}(X_{\rm reg})\neq\emptyset$, hence $\dim\,U=2\dim\,X$.
The inequality $(\dagger)$ for $u$ gives now $\dim\,T_xX=\dim\,X$.
\end{proof}

The last case we consider is that of l.c.i.
toric varieties. If ${\rm char}\,k=0$, then $X$ has rational singularities
by \cite{fulton}, Section 3.5, and it follows by Theorem~\ref{main1}
(see, also, Remark~\ref{rational}) that $X_m$ is irreducible for every $m$.
We give below a direct argument independent of characteristic, which
uses the description due to Nakajima \cite{nakajima}
 of l.c.i. toric varieties. This description
was used by Dais, Haase and Ziegler
 in \cite{dais} to show that all such toric varieties have
crepant resolutions. The main point in our proof is to exhibit a certain
resolution for the ``dual'' toric variety.

\begin{theorem}\label{toric}
If $X$ is an l.c.i. toric variety,
then $X_m$ is irreducible for all $m\geq 1$.
\end{theorem}

\begin{proof}
We use notation and results from \cite{fulton}. Since
all the semigroups we use are saturated, we make no distinction
between the semigroup and the cone it generates.

In general, for two varieties $X$ and $Y$, we have
$(X\times Y)_m\simeq X_m\times Y_m$. Using this, we reduce
immediately to the case when $X=U_{\sigma}$ is affine, where
$\sigma$ is a strongly convex, rational, polyhedral cone
of maximal dimension in $N_{\mathbf R}$, for some lattice
$N$ of rank $n$. Let $S=\sigma^{\vee}\cap M$, where
$M=N^{\vee}$ is the dual lattice.

For every face $\tau\prec\sigma$, we have a corresponding orbit
$O_{\tau}$ of dimension $n-\dim\,\tau$ and a distinguished
point $x_{\tau}\in O_{\tau}$ defined by the semigroup morphism
$x_{\tau}\,:\,S\longrightarrow k$, $x_{\tau}(u)=1$,
if $u\in\tau^{\perp}$ and $x_{\tau}(u)=0$, otherwise.
By Proposition~\ref{irred_crit}, it is enough
to show that if $\{0\}\neq\tau\prec\sigma$, then
$$(\star)\,\dim\pi_m^{-1}(x_{\tau})<mn+\dim\,\tau.$$

We use the inductive description of $S$, due to Nakajima, for the case when
$X$ is locally complete intersection
(see \cite{nakajima}). There are $r\geq 1$ and $s\geq 0$, with
$n=r+s$, such that $S$ can be obtained as follows:
take $S_0=\N^r$, $S\simeq S_s$ and for every $i$, $1\leq i\leq s$,
there is $x\in S_{i-1}\setminus\{0\}$, such that
$$S_i=S_{i-1}+\N e+\N (x-e)\subset S_{i-1}\oplus\ZZ,$$
where $e=(0,1)$.

\smallskip

We show by induction on $s$ that if $T\subset S$ is a face of $S$,
then there is a nonsingular fan $\Delta_s$ 
refining $S$, with rays
$v_1,\ldots,v_{r+s},w_1,\ldots,w_s$ such that:
\item{i)} $\{v_1,\ldots,v_{r+s}\}$ span a cone of $\Delta_s$.
\item{ii)} For every $i\leq s$, $\{w_i,v_{r+j};j\neq i\}$ span a cone
of $\Delta_s$.
\item{iii)} $\{i\leq r+s\vert v_i\in T\}$ has at least
$\dim\,T$ elements.

The assertion is trivial when $s=0$. For $1\leq i\leq s$, let
$\Delta_{i-1}$ be the refinement corresponding to $S_{i-1}$ and 
$T\cap S_{i-1}$. $\Delta_i$ consists of the cones spanned by 
$\{C,e\}$ and $\{C,x-e\}$, where $C$ is a cone in $\Delta_{i-1}$. 
Notice that $\dim(T\cap S_{i-1})\geq\dim(T\cap S_i)-1$, with equality
if and only if $e\in T$ or $x-e\in T$. If we take $v_{r+i}$ to be
$e$ or $x-e$ (we pick the one in $T$, if possible),
and $w_i$ to be the other one, then 
this refinement satisfies the requirements for $S_i$.

\smallskip

In order to prove $(\star)$, notice that $\pi_m^{-1}(x_{\tau})$ consists
of semigroup homomorphisms $\phi\,:\,S\longrightarrow k[t]/(t^{m+1})$,
such that the composition with the projection onto $k$ is $x_{\tau}$.
Let $T=S\cap\tau^{\perp}$ and consider the refinemet $\Delta_s$
constructed above.

Let $a_1,\ldots,a_{r+s}$ be integers such that $1\leq a_i\leq m+1$,
if $v_i\not\in T$, and $a_i=0$, if $v_i\in T$. Take $V_{a_1,\ldots,a_{r+s}}$
to be the locally closed subset of $\pi_m^{-1}(x_{\tau})$ consisting
 of morphisms $\phi$, with 
${\rm ord}\,\phi(v_i)=a_i$. We clearly have $\pi_m^{-1}(x_{\tau})
=\bigcup_{a_1,\ldots,a_{r+s}}V_{a_1,\ldots,a_{r+s}}$.

Fix $i\leq s$. Since $\Delta_s$ is a nonsingular fan, by conditions
$i)$ and $ii)$, there is a relation
$$w_i+v_{r+i}=\sum_{j\leq r+s, j\neq r+i}b_{ij}v_j.$$
This shows that for $\phi\in V_{a_1,\ldots,a_{r+s}}$, the image
of $\phi(w_i)$ in $k[t]/(t^{m+1-a_{r+i}})$ is uniquely determined
by the values of $\phi$ on $v_j$, $j\leq r+s$.
Therefore we get
$$\dim\,V_{a_1,\ldots,a_{r+s}}
\leq\sum_{i=1}^{r+s}(m-a_i+1)-|\{i\leq r+s\vert v_i\in T\}|+\sum_{i=1}^s
a_{r+i}\leq$$
$$(m+1)n-\sum_{i=1}^ra_i-\dim\,T<mn+\dim\tau,$$
since $\dim\,T=n-\dim\tau$, and there is at least one $i\leq r$,
with $v_i\not\in T$ (as $\tau\neq\{0\}$, we have $T\neq S$, and
going downward, we get $T\cap S_0\neq S_0$).
\end{proof}

\vfill\eject

\section*{\bf Appendix}

\bigskip

\begin{center}

{by {\bf David Eisenbud} and {\bf Edward Frenkel}}

\end{center}

\bigskip

Let $G$ be a connected simple algebraic group over an algebraically
closed field $k$ of characteristic zero, and $\g$ its Lie
algebra. Denote by $k[\g]$ the ring of functions on $\g$.  Let $I(\g)$
be the subring of $G$--invariants (equivalently, $\g$--invariants) of
$k[\g]$ under the adjoint action. According to a theorem of B. Kostant
\cite{Ko}, $I(\g) = k[P_1,...,P_\ell]$ for some graded elements
$P_1,...,P_\ell$ of $k[\g]$. Here $\ell = \on{rank} \g$ and $\deg P_i
= d_i$ equals the $i$th exponent of $\g$ plus $1$. By definition, the
nilpotent cone $\NN$ of $\g$ is $\Spe k[\g]/(I(\g)_+)$, where
$(I(\g)_+) = (P_1,...,P_\ell)$ is the ideal generated by the
augmentation ideal $I(\g)_+$ of $I(\g)$.

Kostant \cite{Ko} has proved that $k[\g]$ is free over $I(\g)$. This
important result has numerous applications; in particular, it implies
that the universal enveloping algebra $U(\g)$ is free over its center.

Our goal is to formulate and prove an analogue of Kostant's theorem in
the setting of jet schemes.

Let $G_n$ be the $n$th jet scheme of $G$. This is an algebraic group,
which is isomorphic to a semi-direct product of $G$ and a unipotent
group. The infinite jet scheme $G_\infty$ of $G$ is a proalgebraic
group, which is isomorphic to a semi-direct product of $G$ and a
prounipotent group. The Lie algebra of $G_n$ is the $n$th jet scheme
$\g_n$ of $\g$.  The Lie algebra of $G_\infty$ is the infinite jet
scheme $\g_\infty$ of $\g$. We have: $\g_n = \g[t]/(t^{n+1})$ and
$\g_\infty = \g[[t]]$. The natural projections $\g_\infty \to \g_n$
give us embeddings $\imath_n: k[\g_n] \hookrightarrow k[\g_\infty]$.

Let $d=\dim\g$ and let $x_1,\ldots,x_d$ be a basis of $\g^*$, which we
take as a set of generators of $k[\g]$. We can view them as elements
of $k[\g_n]$ and $k[\g_\infty]$ using the embeddings
$\imath_n$. Define a function $x_i^{(m)}, m\geq 0$, on $\g_n, n\geq m$
and on $\g_\infty$ by the formula
$$
x_i^{(m)}(y(t)) = x_i(\partial_t^m y(t)|_{t=0}), \quad \quad y(t) \in
\g[t]/(t^{n+1}) \on{or} \g[[t]].
$$
The ring $k[\g_n]$ (resp., $k[\g_\infty]$) is the ring of
polynomials in $x_i^{(m)}, 0\leq m\leq n$ (resp., $m\geq 0$). We
introduce a $\Z_{\geq 0}$--gradation on these rings by setting
$\deg x_i^{(m)} = m+1$. Note that each graded component of any of the
above rings is finite-dimensional.

The ring $k[\g_\infty]$ carries a canonical derivation $\pa$ of degree
$1$, defined by the formula $\pa x_i^{(m)} = x_i^{(m+1)}$ (it
corresponds to the vector field $\partial_t$ on $\g[[t]]$). Let
$P_i^{(m)} = \pa^m P_i$. Using the injections $\imath_n$, we view
$P_i^{(m)}$ as elements of $k[\g_n], n\geq m$. We have: $\deg
P_i^{(m)} = d_i+m$. Since the set $\{ P_1,\ldots,P_\ell \}$ is
algebraically independent, it follows that the set $\{
P_1^{(m)},\ldots,P_\ell^{(m)} \}_{m\geq 0}$ is also algebraically
independent.  Because the elements $P_i$ are $G$--invariant, the
element $P_i^{(m)}$ is $G_n$--invariant for $n\geq m$ and
$G_\infty$--invariant.

Denote by $I(\g_n)$ (resp., $I(\g_\infty)$) the ring of
$G_n$--invariants (resp., $G_\infty$--invariants) of $k[\g_n]$ (resp.,
$k[\g_\infty]$) under the adjoint action.

\medskip

\noindent{\bf Proposition A.1.}(\cite{BD}) {\em The ring $I(\g_n)$
(resp., $I(\g_\infty)$) is generated by the algebraically independent
elements $P_1^{(m)},...,P_\ell^{(m)}, 0\leq m\leq n$ (resp., $m\geq
0$). Thus, $I(\g_n)=k[P_1^{(m)},...,P_\ell^{(m)}]_{0\leq m\leq n}$ and
$I(\g_\infty) =$ $k[P_1^{(m)},...,P_\ell^{(m)}]_{m\geq 0}$.}

\medskip

\begin{proof} Let
$\g_{\on{reg}}$ be the smooth open subset of $\g$ consisting of
regular elements. It is known that the morphism $$\chi: \g_{\on{reg}}
\to {\mc P} := \Spe k[P_1,\ldots,P_\ell]$$ is smooth and surjective
(see \cite{Ko}). Therefore the morphism
$$\chi_n: (\g_{\on{reg}})_n \to {\mc P}_n := \Spe
k[P_1^{(m)},...,P_\ell^{(m)}]_{0\leq m\leq n}$$ is also smooth and
surjective.

Consider the map $a: G \times \g_{\on{reg}} \to \g_{\on{reg}}
\underset{{\mc P}}\times \g_{\on{reg}}$ defined by the formula $a(g,x)
= (x,g \cdot x)$. The map $a$ is smooth, and since $G$ acts
transitively along the fibers of $\chi$, it is also surjective. Hence
the corresponding map of jet schemes $a_n: G_n \times
(\g_{\on{reg}})_n \to (\g_{\on{reg}})_n \underset{{\mc P}_n}\times
(\g_{\on{reg}})_n$ is surjective. Given two points $y_1,y_2$ in the
same fiber of $\chi_n$, let $(h,y_1)$ be a point in the (non-empty)
fiber $a_n^{-1}(y_1,y_2)$. Then $y_2 = h \cdot y_1$. Hence $G_n$ acts
transitively along the fibers of the map $\chi_n$.

Since $k[P_1^{(m)},...,P_\ell^{(m)}]_{0\leq m\leq n}$ is normal, this
implies that the ring of $G_n$--invariant functions on
$(\g_{\on{reg}})_n$ equals $k[P_1^{(m)},...,P_\ell^{(m)}]_{0\leq m\leq
n}$. Because $\g_{\on{reg}}$ is a smooth open subset of $\g$, we
obtain that $(\g_{\on{reg}})_n$ is dense in $\g_n$, and so any
$G_n$--invariant function on $\g_n$ is determined by its restriction
to $(\g_{\on{reg}})_n$. This proves the proposition in the case of the
finite jet schemes. The same argument works in the case the infinite
jet scheme as well.
\end{proof}

Let $I(\g_{n})_+$ be the augmentation ideal of the graded ring
$I(\g_n)$. By Proposition A.1, the ideal $(I(\g_{n})_+)$ in $k[\g_n]$
generated by $I(\g_n)_+$ equals $(P_1^{(m)},...,P_\ell^{(m)})_{0\leq
m\leq n}$. Hence we obtain that the $n$th jet scheme $\NN_n$ of the
nilpotent cone $\NN$ is $\Spe k[\g_n]/(I(\g_{n})_+)$.
Likewise, $\NN_\infty =$ $\Spe k[\g_{\infty}]/(I(\g_{\infty})_+)$.

According to Kostant's results \cite{Ko} (see \cite{BL}, \cite{CG} for a
review), the nilpotent cone is a complete intersection, which is
irreducible and reduced. Moreover, it is proved in \cite{He} that it
has rational (hence canonical) singularities. Therefore we obtain from
\thmref{main1} and \propref{reduced}:

\medskip

\noindent{\bf Theorem A.2.}    \label{one}
{\em $\NN_n$ is irreducible, reduced and a complete intersection.}

\medskip

\noindent{\bf Corollary A.3.}    \label{dom}
{\em The natural map $k[\NN_n] \to k[\NN_{n+1}]$ is an embedding.}

\medskip

\begin{proof}
Let $Y$ be the open dense $G$--orbit of regular elements in $\NN$. By
Theorem A.2, $\NN_n$ is irreducible. Hence $Y_n$ is dense in
$\NN_n$. Since $Y$ is smooth, $Y_{n+1} \to Y_n$ is
surjective. Therefore the map $\NN_{n+1} \to \NN_n$ is dominant.
\end{proof}

Now we can formulate an analogue of the Kostant freeness theorem.

\medskip

\noindent{\bf Theorem A.4.}    \label{two}
(1) {\em $k[\g_n]$ is free over $I(\g_n)$.}

(2) {\em $k[\g_\infty]$ is free over $I(\g_\infty)$.}

\medskip

\begin{proof}
Since $\NN_n$ is a complete intersection, it is
Cohen-Macaulay. Therefore $k[\g_n]$ is a flat module over $I(\g_n) =
k[P_1^{(m)},...,P_\ell^{(m)}]_{0\leq m\leq n}$. Since $I(\g_n)$ is
$\Z_+$--graded with finite-dimensional homogeneous components,
flatness implies that $k[\g_n]$ is free over $I(\g_n)$ (see, e.g.,
Ex. 18.18 of \cite{eisenbud}). This proves part (1).

Since $k[\g_{n}]$ is free over $I(\g_{n})$, and both rings are
$\Z_+$--graded with finite-dimen\-sional homogeneous components, we
obtain that any graded lifting of a $k$--basis of $k[\NN_n]$ to
$k[\g_n]$ is an $I(\g_n)$--basis of $k[\g_n]$. Conversely, the image
of any $I(\g_n)$--basis of $k[\g_n]$ in $k[\NN_n]$ is a $k$--basis of
$k[\NN_n]$.

Now choose any graded basis $S_n$ of $k[\g_n]$ over $I(\g_n)$.  Then
the image $S'_n$ of $S_n$ in $k[\NN_n]$ is a $k$--basis of
$k[\NN_n]$. According to Corollary A.3, the image of $S'_n \subset
k[\NN_n]$ in $k[\NN_{n+1}]$ can be extended to a $k$--basis of
$k[\NN_{n+1}]$. Hence the image of $S_n$ in $k[\g_{n+1}]$ (we denote
it also by $S_n$) can be extended to an $I(\g_{n+1})$--basis $S_{n+1}$
of $k[\g_{n+1}]$.

Thus, we obtain a directed system $S_n, n\geq 0$, of sets of graded
elements of $k[\g_{\infty}]$ and embeddings $S_n \hookrightarrow S_m,
\forall m\geq n$. Let $S$ be the union of all of the sets $S_n, n\geq
0$. Note that by Proposition A.1, $I(\g_\infty)$ is the inductive
limit of the directed system $I(\g_n), n\geq 0$. We claim that $S$ is
a basis of $k[\g_\infty]$ over $I(\g_\infty)$.

Indeed, consider the multiplication map $m: \on{span}(S) \otimes
I(\g_\infty) \to k[\g_\infty]$. This map is surjective: take any
homogeneous element $A$ of $\g_\infty$ of degree $a$. Then it already
belongs to $k[\g_a]$. Hence by construction $A$ lies in the image of
$m$. The map $m$ is also injective. Indeed, suppose there is an
element $B$ in the kernel of $m$. Without loss of generality we can
assume that $B$ is homogeneous of degree $b$. But then $B$ belongs
to the kernel of the map $\on{span}(S_b) \otimes I(\g_b) \to k[\g_b]$,
and so $B=0$ by construction. This completes the proof.
\end{proof}

\noindent{\bf Remark A.5.}
Part (1) of Theorem A.4 has been previously proved by Geoffriau
\cite{Ge} in the case when $n=1$ (for arbitrary $\g$) and in the case
when $\g=\sw_2$ (for arbitrary $n$). We thank M. Duflo for bringing
the paper \cite{Ge} to our attention.

\medskip

\noindent{\bf Remark A.6.} Bernstein and Lunts \cite{BL} have given an
alternative proof of the original Kostant freeness theorem, using the
isomorphism $I(\g) \simeq k[\hh]^W$, and the Chevalley theorem which
states that $k[\hh]$ is a free over $k[\hh]^W$ (here $\hh$ denotes the
Cartan subalgebra of $\g$, and $W$ denotes the Weyl group of $\g$).
In our case this approach cannot be applied: although there is a
natural embedding of $I(\g_n)$ to $k[\hh_n]$, the ring $k[\hh_n]$ is
not free over the image of $I(\g_n)$. In fact, the latter equals
$k[(\hh/W)_n]$, which is strictly smaller than $k[\hh_n]^W$.

\medskip

Theorem A.2 also implies the following result:

\medskip

\noindent{\bf Proposition A.7.}    \label{three}
{\em The space of $G_n$--invariants (resp., $G_\infty$--invariants) of
$k[\NN_n]$ (resp., $k[\NN_\infty]$) consists of constants.}

\begin{proof}
The jet scheme $Y_n$ (resp, $Y_\infty$) of the open dense $G$--orbit
$Y$ of $\NN$ is an orbit of the group $G_n$ (resp., $G_\infty$) in
$\NN_n$ (resp., $\NN_\infty$). Therefore any invariant function on it
is a constant. But according to the proof of Corollary A.3, it is a
dense subvariety in $\NN_n$ (resp., $\NN_\infty$). Hence any invariant
function on $\NN_n$ (resp., $\NN_\infty$) is a constant.
\end{proof}

\medskip

\noindent{\bf Remark A.8.}  According to Theorem A.4, the natural
morphisms $\g_\infty \to \Spe I(\g_\infty)$ and $\g_n\to\Spe I(\g_n)$
are flat. Drinfeld has suggested that these morphisms may be viewed as
local counterparts of the Hitchin morphism. More precisely, let $X$ be
a smooth projective curve over $\C$, and $\on{Bun}_G$ the moduli stack
of $G$--bundles on $X$. The cotangent space $T^*_{\mc F} \on{Bun}_G$
to $\on{Bun}_G$ at ${\mc F} \in \on{Bun}_G$ is isomorphic to
$H^0(X,\g_{\mc F} \otimes \Omega)$. Here $\g_{\mc F} = {\mc F}
\underset{G}\times \g$, and we identify $\g \simeq \g^*$ using the
invariant inner product on $\g$. The Hitchin morphism $$T^* \on{Bun}_G
\to {\mc H} = \oplus_{i=1}^\ell H^0(X,\Omega^{\otimes d_i})$$ sends
$({\mc F},\omega \in H^0(X,\g_{\mc F} \otimes \Omega))$ to
$(P_1(\omega),\ldots,P_\ell(\omega)) \in {\mc H}$.

Let $x$ be a point of $X$, and $\wh{\mc O}_x$ the completion of the
local ring at $x$. Denote by $D_x$ the formal disc at $x$, $D_x = \Spe
\wh{\mc O}_x$. For each ${\mc F} \in \on{Bun}_G$ we have a local
analogue of the Hitchin map, $$h_x^{\mc F}: H^0(D_x,\g_{\mc F} \otimes
\Omega) \to {\mc H}_x = \oplus_{i=1}^\ell H^0(D_x,\Omega^{\otimes
d_i}).$$ If we trivialize ${\mc F}|_{D_x}$ and choose a formal
coordinate $t$ at $x$, the map $h_x^{\mc F}$ becomes our map
$\g_\infty \to \Spe I(\g_\infty)$. Actually, for varying ${\mc F}$ and
$x \in X$, the spaces $H^0(D_x,\g_{\mc F} \otimes \Omega)$ and ${\mc
H}_x$ can be glued together into schemes over $X \times \on{Bun}_G$
equipped with a flat connection along $X$, and the maps $h_x^{\mc F}$
can be glued into a morphism between these schemes preserving
connections. The Hitchin morphism then appears as the corresponding
map of the schemes of horizontal sections.

The flatness of the Hitchin morphism has been proved by Hitchin
\cite{Hi} and Faltings \cite{Fa}. Drinfeld has derived the flatness of
the morphism $\g_n\to\Spe I(\g_n)$ from the flatness of the Hitchin
morphism (private communication). He asked whether one can find a
purely ``local'' argument proving this fact. The above proof answers
this question.


\providecommand{\bysame}{\leavevmode\hbox to3em{\hrulefill}\thinspace}


\bigskip
\vbox{\noindent Author Addresses:}
\smallskip
\noindent{David Eisenbud}\par
\noindent{Department of Mathematics, University of California, Berkeley,
CA, 94720; eisenbud@math.berkeley.edu}\par
\smallskip
\noindent{Edward Frenkel}\par
\noindent{Department of Mathematics, University of California, Berkeley,
CA, 94720; frenkel@math.berkeley.edu}\par
\smallskip
\noindent{Mircea Musta\c{t}\v{a}}\par
\noindent{Department of Mathematics, University of California, Berkeley,
CA, 94720 and Institute of Mathematics of the Romanian Academy; 
mustata@math.berkeley.edu}\par

\end{document}